\documentclass{article}

\usepackage{amssymb}
\usepackage{amsmath}

\usepackage[dvips]{graphicx}

\begin{document}


\begin{center}
\textbf{A remark on $n-$ Jordan homomorphisms }
\end{center}
 
\noindent \textbf{}

\begin{center}
M. El Azhari 
\end{center}

\noindent \textbf{}

\noindent \textbf{}

\noindent \textbf{Abstract.} Let $ A $ and $ B $ be commutative algebras and $ n\geqslant 2 $ an integer. Then each $n-$ Jordan homomorphism $ h:A\rightarrow B $ is an $n-$ homomorphism.

\noindent \textbf{}


\noindent \textbf{} 


\noindent \textbf{}

\noindent \textbf{}

\noindent \textbf{1. Preliminaries }

\noindent \textbf{}                                                                        

\noindent\textbf{} Let $ A, B $  be two algebras and $ n\geqslant 2 $ an integer. A linear map $ h:A\rightarrow B $ is called an $n-$ Jordan homomorphism if $ h(a^{n})=h(a)^{n} $ for all  $ a\in A. $ Also, a linear map $ h:A\rightarrow B $  is called an $n-$ homomorphism if $ h(\prod_{i=1}^{n}a_{i})=\prod_{i=1}^{n} h(a_{i}) $ for all  $ a_{1},...,a_{n}\in A. $ In the usual sense, a $ \:2-$ Jordan homomorphism is a Jordan homomorphism and a $ \:2-$ homomorphism is a homomorphism.\\
It is obvious that $ n- $  homomorphisms are $ n- $ Jordan homomorphisms. Conversely, under certain conditions, $ n- $ Jordan homomorphisms are $ n- $  homomorphisms. For example, if $ A, B $ are commutative algebras and $ h:A\rightarrow B $ is a Jordan homomorphism, then $ h $ is a homomorphism: Let $ x, y\in A, $ since $ h((x+y)^{2})=h(x+y)^{2}, $ a simple calculation shows that $ h(xy)=h(x)h(y). $ Also, it was shown in [3] that if $ n\in\{3,4\}, A, B $ are commutative algebras and $ h:A\rightarrow B $ is an $n-$ Jordan homomorphism, then $ h $ is  an $n-$ homomorphism. This result was also proved for $ n=5 $ in [4].                                                                                                                                                                                                                                                                                                                                                                                                                                                                                                                                                                                                                                                                                                                                                                                                                                                                                                                                                                                                                                                                                                                                                                                                                                                                                                                                                                                                                                                                                                                                                                            

\noindent \textbf{}

\noindent \textbf{2. Results} 

\noindent \textbf{}

\noindent \textbf{} Cheshmavar et al. [2] claimed the following result (see also [1], [5] and [6]):

\noindent \textbf{}

\noindent \textbf{ } Theorem 2.3 of [2]. Let  $ A $ and $ B $ be commutative algebras, $ n\geqslant 3, $ and let  $ h:A\rightarrow B $ be an $ n- $ Jordan homomorphism. Then $ h $ is an $ n- $ homomorphism. 

\noindent \textbf{}

\noindent \textbf{} The proof of this theorem is essentially based on a formula which is not correct: Let  $ h:A\rightarrow B $ be a linear map, $ n\geqslant 4, $  and consider the maps\\
$ \psi: A^{n}\rightarrow B,\:\psi(x_{1},...,x_{n})= h(\prod_{i=1}^{n}x_{i})-\prod_{i=1}^{n} h(x_{i})$\\
and for $ 2\leqslant k\leqslant n, $\\
$ \varphi_{k-1}: A^{k}\rightarrow B,\:\varphi_{k-1}(x_{1},...,x_{k})=h((x_{1}+...+x_{k})^{n})-h(x_{1}+...+x_{k})^{n}.$\\
The authors [2] claimed that\\
\\
$ \varphi_{n-1}(x_{1},...,x_{n})=\sum_{i,j=1,i<j}^{n}\varphi_{1}(x_{i},x_{j})+\sum_{i,j,k=1,i<j<k}^{n}\varphi_{2}(x_{i},x_{j},x_{k})\\
+...+\sum_{i_{1},...,i_{n-1}=1,i_{1}<...<i_{n-1}}^{n}\varphi_{n-2}(x_{i_{1}},...,x_{i_{n-1}})+n!\:\psi(x_{1},...,x_{n}). $\\
\\
This formula is false since, for example, the expression containing exactly the coordinates $ x_{1},x_{2} $ in the left part will be appearing many times in the right part ( in $ \varphi_{1}(x_{1},x_{2}),\varphi_{2}(x_{1},x_{2},x_{3}),\varphi_{2}(x_{1},x_{2},x_{4}),...$ ).

\noindent \textbf{}

\noindent \textbf{} We give here a correct proof of [2, Theorem 2.3].

\noindent \textbf{}

\noindent \textbf{Theorem 2.1.} Let $ A, B $ be two algebras, $ n\geqslant 2, $ and let  $ h:A\rightarrow B $ be an $ n- $ Jordan homomorphism. Then $ \sum_{\sigma\in S_{n}} (h(\prod_{i=1}^{n}x_{\sigma(i)})-\prod_{i=1}^{n} h(x_{\sigma(i)}))=0 $ for all $ x_{1},...,x_{n}\in A, $ where $ S_{n} $ is the set of permutations on $ \{1,...,n\}. $ Furthermore, if $ A $ and $ B $ are commutative, then $ h $
is an $ n- $ homomorphism.

\noindent \textbf{}

\noindent \textbf{Proof.} Case $ n=2:$ Let $ h:A\rightarrow B $ be a Jordan homomorphism and $ x_{1}, x_{2} \in A.$ Since $ h((x_{1}+x_{2})^{2}) = h(x_{1}+x_{2})^{2},$ a simple calculation shows that $h(x_{1}x_{2})-h(x_{1})h(x_{2})+h(x_{2}x_{1})-h(x_{2})h(x_{1})=0,$ 
then  \\
$ \sum_{\sigma\in S_{2}} (h(x_{\sigma(1)}x_{\sigma(2)})-h(x_{\sigma(1)})h(x_{\sigma(2)}))=0. $ \\
If $ A $ and $ B $ are commutative, we deduce that $ h $ is a homomorphism.\\
\\
Case $ n\geqslant 3: $ Let $ h:A\rightarrow B $ be a linear map and $ x_{1},...,x_{n}\in A. $ We will write $ h((x_{1}+...+x_{n})^{n})-h(x_{1}+...+x_{n})^{n} $ as the sum of expressions containing respectively $ 1,2,...,n $ coordinates of $ (x_{1},...,x_{n}). $ Let $ 1\leqslant k\leqslant n $ and $ i_{1},...,i_{k}\in \{1,...,n\},\: i_{1}<...<i_{k},$ we denote by $ \psi_{k}^{(n)}(x_{i_{1}},...,x_{i_{k}})$ the expression in $ h((x_{1}+...+x_{n})^{n})-h(x_{1}+...+x_{n})^{n} $ containing exactly the $ k $
coordinates $ x_{i_{1}},...,x_{i_{k}}. $ Then\\
$ h((x_{1}+...+x_{n})^{n})-h(x_{1}+...+x_{n})^{n}=\sum_{i=1}^{n}\psi_{1}^{(n)}(x_{i})+\sum_{i,j=1,i<j}^{n}\psi_{2}^{(n)}(x_{i},x_{j})\\
+...+\sum_{i_{1},...,i_{n-1}=1,i_{1}<...<i_{n-1}}^{n}\psi_{n-1}^{(n)}(x_{i_{1}},...,x_{i_{n-1}})+\psi_{n}^{(n)}(x_{1},...,x_{n}). $\\
\\
For $ 1\leqslant i\leqslant n,\:\psi_{1}^{(n)} (x_{i})=h(x_{i}^{n})-h(x_{i})^{n}.$ \\
\\
For $ i,j\in\{1,...,n\},\:i<j,\\
\psi_{2}^{(n)}(x_{i},x_{j})=h((x_{i}+x_{j})^{n}) - h(x_{i}+x_{j})^{n}-(h(x_{i}^{n}) - h(x_{i})^{n})-(h(x_{j}^{n}) - h(x_{j})^{n})\\
=h((x_{i}+x_{j})^{n}) - h(x_{i}+x_{j})^{n}-\psi_{1}^{(n)}(x_{i})-\psi_{1}^{(n)}(x_{j}).$\\
\\
For $ 3\leqslant k\leqslant n,\:i_{1},...,i_{k}\in\{1,...,n\}, i_{1}<...<i_{k}, \\
\psi_{k}^{(n)}(x_{i_{1}},...,x_{i_{k}})=  h((x_{i_{1}}+...+x_{i_{k}})^{n})-h(x_{i_{1}}+...+x_{i_{k}})^{n}\\
-\sum_{m=1}^{k}\psi_{1}^{(n)}(x_{i_{m}})-\sum_{j_{1},j_{2}\in\{ i_{1},...,i_{k}\},j_{1}<j_{2}}\psi_{2}^{(n)}(x_{j_{1}},x_{j_{2}})\\
-...-\sum_{j_{1},...,j_{k-1}\in\{ i_{1},...,i_{k}\},j_{1}<...<j_{k-1}}\psi_{k-1}^{(n)}(x_{j_{1}},...,x_{j_{k-1}}). $\\
\\
If  $ h:A\rightarrow B $ is an $ n- $ Jordan homomorphism, then $ \psi_{1}^{(n)}(x_{i})=0 $ for\\
 $ 1\leqslant i\leqslant n, $ also $ \psi_{2}^{(n)}(x_{i},x_{j})=0 $ for $  i,j\in\{1,...,n\},\:i<j, $ and by induction
$ \psi_{n}^{(n)}(x_{1},...,x_{n})= \sum_{\sigma\in S_{n}} (h(\prod_{i=1}^{n}x_{\sigma(i)})-\prod_{i=1}^{n} h(x_{\sigma(i)}))=0. $\\
\\
If $ A $ and $ B $ are commutative,\\
\\
$  \sum_{\sigma\in S_{n}} (h(\prod_{i=1}^{n}x_{\sigma(i)})-\prod_{i=1}^{n} h(x_{\sigma(i)}))=
n!\:(h(\prod_{i=1}^{n}x_{i})-\prod_{i=1}^{n} h(x_{i}))=0,$\\
\\
hence $ h $ is an $ n- $ homomorphism.

\noindent \textbf{}

\begin{center}
REFERENCES
\end{center}
 
\noindent \textbf{} 

\noindent \textbf{} [1] A. Bodaghi and H. Inceboz, $ n- $ Jordan homomorphisms on commutative algebras, Acta Math. Univ. Comenianae, 87(1)(2018), 141-146.

\noindent \textbf{} [2] J. Cheshmavar, S. K. Hosseini, and C. Park, Some resuls on $ n- $ Jordan homomorphisms, Bull. Korean Math. Soc. 57(1)(2020), 31-35.

\noindent \textbf{} [3] M. Eshaghi Gordji, $ n- $ Jordan homomorphisms, Bull. Aust. Math. Soc. 80(1)(2009), 159-164.

\noindent \textbf{} [4] M. Eshaghi Gordgi, T. Karimi, and S. Kaboli Gharetapeh, Approximately $ n- $ Jordan homomorphisms on Banach algebras, J. Inequal. Appl. 2009(2009), Article ID 870843, 8 pages.

\noindent \textbf{} [5] E. Gselmann, On approximate $ n- $ Jordan homomorphisms, Annales Mathematicae Silesianae, 28(2014), 47-58.

\noindent \textbf{} [6] Yang-Hi Lee, Stability of $ n- $ Jordan homomorphisms from a normed algebra to a Banach algebra, Abstract and Applied Analysis, 2013(2013), Article ID 691025, 5 pages.

\noindent \textbf{} 
 
\noindent \textbf{} 

\noindent \textbf{} Ecole Normale Sup\'{e}rieure

\noindent \textbf{} Avenue Oued Akreuch

\noindent \textbf{} Takaddoum, BP 5118, Rabat

\noindent \textbf{} Morocco
 
\noindent \textbf{} 

\noindent \textbf{} E-mail:  mohammed.elazhari@yahoo.fr

\end{document}